\documentclass[11pt,oneside]{amsart}
\usepackage{amsmath,amsthm,amssymb,amscd,fancybox,esint}
\usepackage[a4paper,left=3cm,right=3cm,top=3cm,bottom=3cm]{geometry}
\usepackage{cite}
\usepackage{mathtools}
\renewcommand{\epsilon}{\varepsilon}
\renewcommand{\phi}{\varphi}

 \newcommand{\bZ}{\mathbb{Z}}
\newcommand{\bR}{\mathbb{R}}

\newcommand{\ddS}{\mathcal{S}'(\bR^n)}  

\newtheorem{theorem}{Theorem}[section]

\newtheorem{proposition}[theorem]{Proposition}

\theoremstyle{definition}
\newtheorem{definition}[theorem]{Definition}

\theoremstyle{remark}
\newtheorem{remark}[theorem]{Remark}
\numberwithin{equation}{section}

\def\bZ{{\mathbb Z}}

\def\bC{{\mathbb C}}
\def\bR{{\mathbb R}}

\def\cQ{\mathcal{B}}

\def\cQ{\mathcal{Q}}

\newcommand{\vf}{\mathbf{f}}
\newcommand{\vg}{\mathbf{g}}
\newcommand{\vx}{\mathbf{x}}
\newcommand{\vw}{\mathbf{w}}

\newcommand{\C}{\mathbb{C}}
\newcommand{\cD}{{\mathbb{R}^n}}

\def\cprime{$'$}
\newcommand{\esssup}{\mathop{\rm ess{\,}sup}}

\providecommand{\norm}[2][]{#1\lVert#2#1\rVert}
\providecommand{\abs}[2][]{#1\lvert#2#1\rvert}
\allowdisplaybreaks
\begin{document}
\title[]{ Bandlimited multipliers on matrix-weighted $L^p$-spaces }
\author[M.\ Nielsen]{Morten Nielsen}
\address{Department of Mathematical Sciences\\ Aalborg
  University\\ Thomas Manns Vej 23\\ DK-9220 Aalborg East\\ Denmark}
\email{mnielsen@math.aau.dk}
\begin{abstract}
 We extend a classical result by Triebel on boundedness of  bandlimited multipliers on $L^p(\bR^n)$, $0<p\leq 1$, to  a vector-valued and matrix-weighted setting with boundedness of the bandlimited multipliers obtained on $L^p(W)$, $0<p\leq 1$, for  matrix-weights $W:\bR^n\rightarrow \C^{N\times N}$ that satisfy a matrix Muckenhoupt $A_p$-condition.  
\end{abstract}
\maketitle

\section{Introduction}
An $N\times N$ matrix weight on $\bR^n$ is a locally integrable and almost everywhere positive definite matrix function $W\colon \cD\rightarrow \C^{N\times N}$. 
The matrix-weighted $L^p$-space $L^p(W)$, $0< p<\infty$, is defined
for any matrix-weight $W\colon
\cD\to\C^{N\times N}$  as  the family of
measurable functions $\vf\colon \cD\to\C^N$ satisfying 
\begin{equation}\label{eq:lp}
  \norm{\vf}_{L^p(W)}:=\Biggl(\int_{\cD}\abs{W^{1/p}(x)\vf(x)}^p\,dx\Biggr)^{1/p}<\infty.
\end{equation}
Using the standard identification of (vector-)functions that differ on a set of measure zero, one can verify that $L^p(W)$ is  a (quasi-)Banach space. 

The matrix-weighted $L^p$-spaces have attracted a great deal of  attention recently (see, e.g., \cite{MR4454483,MR3687948,Rou03a,MR1857041}) due to the fact that the setup  generates a number of interesting mathematical questions for vector valued functions that  naturally connect to classical results on Muckenhoupt weights in harmonic analysis. A highlight 
in the matrix-weighted case is the formulation of a suitable matrix $A_p$ condition by Nazarov, Treil and Volberg that completely characterizes boundedness of the Riesz-transform(s) on $L^p(W)$ for $1<p<\infty$, see \cite{Vol97a,TreVol97a}, see also \cite{Gol03a}.

%
In the present paper we study Fourier multipliers on $L^p(W)$ with a focus on the case $0<p\leq 1$. As is well-known, a scalar Fourier multiplier is a function $\phi\in L^\infty(\bR^n)$ that induces a corresponding bounded multiplier operator
$$\phi(D)f:=\mathcal{F}^{-1}(\phi\mathcal{F}f),\qquad f\in L^2(\bR^n),$$
where $\mathcal{F}$ denotes the Fourier transform on $L^2(\bR^n)$, where we used the normalisation specified in Eq. \eqref{def:f} below. 

 In case, $\mathcal{F}^{-1}\phi\in L^1(\bR^n)$, Young's inequality provides an easy extension of the Fourier multiplier to $L^p(\bR^n)$, $1\leq p\leq \infty$, through the estimate
\begin{equation}\label{eq:young}
\|\phi(D)f\|_{L^p(\bR^n)}\leq C\|\mathcal{F}^{-1}\phi\|_{L^1(\bR^n)}\|f\|_{L^p(\bR^n)}.
\end{equation}
Let us now lift the multiplier to the vector-setting. Given a suitably nice  $\vf\colon \cD\to\C^N$, we let  $\phi(D)$ act coordinate-wise on $\vf=(f_1,\ldots,f_N)^T$, i.e.,
$\phi(D)\vf:=(\phi(D)f_1,\ldots, \phi(D)f_N)^T$. A much more challenging question is then whether $\phi(D)$ extends to a bounded operator on $L^p(W)$ for a given matrix-weight $W:\bR^n\rightarrow\bC^{N\times N}$? This problem has been studied in detail in the case $1\leq p<\infty$ by various authors. In case $1<p<\infty$, the results on singular integrals obtain by Goldberg in \cite{Gol03a} provides boundedness $\phi(D):L^p(W)\rightarrow L^p(W)$ provided that $\mathcal{F}^{-1}\phi$ satisfies a mild decay condition, and, more importantly, also provided that the matrix-weight $W$ satisfies the so-called matrix Muckenhoupt $A_p$-condition, see Definition \ref{def:mu} below. Frazier and Roudenko extended this result in \cite{MR4263690} to obtain boundedness in the end-point case $\phi(D):L^1(W)\rightarrow L^1(W)$, provided $W$ satisfies a matrix Muckenhoupt $A_1$-condition. This leaves the range $0<p<1$, which we will study in the present article.

Even in the scalar case, it is known that there are no unrestricted extensions of Young's estimate \eqref{eq:young} to the range $0<p<1$, so let us discuss a  framework that makes the case $0<p<1$ manageable. Here the notion of bandlimited functions will be central. Triebel showed \cite[Theorem 1.5.1]{Triebel1983} that in case  $0<p\leq 1$,  $\phi:B(0,1)\rightarrow\bC$ is compactly supported, and $f$ is a bandlimited tempered distribution in the sense that the frequency support of $f$ satisfies $\text{supp}(\hat{f})\subseteq \{x\in\bR^n: |x|<1\}$, there exists $C$ independent of $f$ such that
$$\|\phi(D)f\|_{L^p(\bR^n)}\leq C\|f\|_{L^p(\bR^n)},$$
 provided  $\mathcal{F}^{-1}\phi$ satisfies a mild decay condition.

We will extend Triebel's result on bandlimited multipliers to the matrix-weighted case in Section \ref{sec:main} for weights $W$ that satisfy a matrix $A_p$-condition.
Finally, in Section \ref{sec:exa}, we provide a specific motivation for studying vector-valued Fourier multipliers in the case $0<p\leq 1$ by considering an application of the multiplier result to the study of matrix-weighted smoothness spaces.

\section{Muckenhoupt weights and Fourier multipliers}\label{sec:main}

The matrix $A_p$-conditions will be of fundamental importance in order to derive our main multiplier result, so let us first define these conditions.  We let $\cQ$ denote the collection of all cubes  $\{Q(z,r)\}_{z\in\bR^n,r>0}$ in $\bR^n$, where $Q(z,r):=z+r[-1/2,1/2)^n$. 
\begin{definition}\label{def:mu}
Let $W\colon \cD\rightarrow \C^{N\times N}$ be a matrix weight. We say that $W$ satisfies  the matrix Muckenhoupt  $A_p$ condition for  $1<p<\infty$ provided
\begin{equation}\label{eq:Roudenko}
 [W]_{{\mathbf{A}_p(\cD)}}:=\sup_{Q\in \cQ} \int_Q\left( \int_Q \big\|W^{1/p}(x)W^{-1/p}(t)\big\|^{p'} \frac{dt}{|Q|}\right)^{p/p'} \frac{dx}{|Q|}<+\infty.
\end{equation}
In case $0<p\leq 1$,  $W$ is said to satisfy  the matrix Muckenhoupt  $A_p$ condition provided 
\begin{equation}\label{eq:mA1}
	[W]_{\mathbf{A}_p(\cD)}:=\sup_{Q\in\cQ} \esssup_{y\in Q} \frac{1}{|Q|}\int_Q \|W(t)^{1/p}W^{-1/p}(y)\|^p\,dt<+\infty.
\end{equation}
The norm $\|\cdot\|$ appearing in the integrals is any matrix norm on the $N\times N$ matrices. In case either \eqref{eq:Roudenko} or \eqref{eq:mA1} applies, we write $W\in \mathbf{A}_p(\cD)$.

\end{definition}
\begin{remark}
The matrix $A_p$-condition for $1<p<\infty$ was introduced in \cite{NazTre96a,TreVol97a,Vol97a} using the notion of dual norms. The condition given in \eqref{eq:Roudenko} was first studied by  Roudenko in \cite{Rou03a},
  where the condition is also proven to  be equivalent to the original matrix $A_p$-condition.    For $0<p\leq 1$,  Frazier and Roudenko introduced the condition \eqref{eq:mA1} 
in  \cite{Frazier:2004ub}. Also, in the scalar case $N=1$, one can verify for $1\leq p<\infty$ that the conditions in Definition \ref{def:mu} are equivalent to the well known scalar $A_p$ conditions.
\end{remark}
\begin{remark}\label{rem:inv}
Since the family of cubes $\cQ$ in $\bR^n$ is clearly invariant under dilations $\vx\rightarrow R\vx$ for $R>0$, it easily follows from \eqref{eq:Roudenko} and \eqref{eq:mA1} that $\mathbf{A}_p(\cD)$ is invariant under such dilations. In fact, for any $W\in \mathbf{A}_p(\cD)$, we have 
$$[W(R\cdot)]_{\mathbf{A}_p(\cD)}=[W]_{\mathbf{A}_p(\cD)}.$$
\end{remark}

 Let us suppose that $W\in \mathbf{A}_p(\cD)$. For a fixed vector $\vx\in\bC^N$, consider the scalar weight $w_{\vx}(t):=\|W^{1/p}(t)\vx\|^p$. For $1<p<\infty$, it is know that $w_{\vx}$ is a scalar $A_p$-weight with $A_p$ constant depending only on $[W]_{\mathbf{A}_p(\cD)}$, see \cite[Corollary 2.2]{Gol03a}, and for 
$0<p\leq 1$, it was shown in  \cite[Lemma 2.1]{Frazier:2004ub} that $w_{\vx}$ is in scalar $A_1$ with an  $A_1$ constant that only depends on $[W]_{\mathbf{A}_p(\cD)}$. In both cases, we may conclude, see, e.g., \cite[Proposition 9.1.5.]{Gra14b}, that $w_{\vx}$ induces a doubling measure in the sense that there exists a constant $C$ independent of $\vx$ such that for any $Q(z,r)\in\cQ$,
\begin{equation}\label{eq:dou}
\int_{Q(z,2 r)} w_{\vx}(t)\,dt\leq C\int_{Q(z,r)} w_{\vx}(t)\,dt.
\end{equation}
The doubling exponent $\beta$ is defined by letting $2^\beta=C$, with $C$ the smallest value of $C$ satisfying \eqref{eq:dou}. It is know that $\beta\geq n$, see, e.g., \cite[Proposition 2.10]{MR3330610}.
We observe that for the unit cubes $Q_k:=Q(k,1)$, $k\in \bZ^n$, the doubling condition ensures that there exists a constant $c>0$, independent of $\vx$, such that for $k,\ell\in\bZ^n$,
$$\int_{Q_k} w_{\vx}(t)\,dt\leq c(1+|k-\ell|)^\beta \int_{Q_\ell} w_{\vx}(t)\,dt,$$
with $\beta$ the doubling exponent of $w_{\vx}$, due to the fact that $Q(k,1)\subseteq Q\big(\ell,2\sqrt{n}(1+|k-\ell|)\big)$.


Let us specify our chosen normalisation of the Fourier transform. For  $f\in L^1(\bR^n)$, we let 
\begin{equation}\label{def:f}\mathcal{F}(f)(\xi):=(2\pi)^{-n/2}\int_{\bR^n}
f(x)e^{- i x\cdot\xi}\,dx,\qquad \xi\in\bR^n,
\end{equation}
denote the Fourier transform, and we use the standard notation $\hat{f}(\xi)=\mathcal{F}(f)(\xi)$. With this normalisation, the Fourier transform extends to a unitary transform on $L^2(\bR^n)$ and we denote the inverse Fourier transform by $\mathcal{F}^{-1}$.

For $R>0$, we define the following class of vector-valued functions with band-limited coordinate functions, where $B(a,r)$ denotes the (open) Euclidean ball in $\bR^n$ centered at $a$ with radius $r$, and $ \ddS$ denotes the tempered distributions defined on $\bR^n$,
\begin{equation}
E_R:=\{\vf:\bR^n\rightarrow \bC^N: f_i\in \ddS \text{ and } \text{supp}(\hat{f}_i)\subseteq B(0,R), i=1,\ldots,N\}.
\end{equation}
We are now ready to state and prove our main result.
\begin{proposition}\label{prop:main}
Let $W\in\mathbf{A}_p(\bR^n)$ for some $0<p\leq 1$ and let $\beta>0$ be the doubling exponent from Eq.\ \eqref{eq:dou} associated with $W$.  Suppose there is a constant $K$ such that the compactly supported function $\phi:B(0,R)\rightarrow\bC$ satisfies $$|\mathcal{F}^{-1}(\phi)(x)|\leq KR^{n} (1+R|x|)^{-M},\qquad x\in\bR^n,$$ for some $M>(n+\beta)/p$. Then there exists a finite constant $C:=C([W]_{\mathbf{A}_p(\bR^n)},K,p)$ such that the Fourier multiplier
$$\phi(D)f:=\mathcal{F}^{-1}[\phi\cdot \mathcal{F}(f)],$$ defined for $f\in \ddS $ with $\text{supp}(\hat{f})\subseteq B(0,R)$, satisfies
$$\|\phi(D)\vf\|_{L^p(W)}\leq C \|\vf\|_{L^p(W)}$$
for all $\vf\in E_R\cap L^p(W)$.
\end{proposition}
\begin{proof}
Let us first consider the special case $R=1$. For any $\vf\in E_1$, we notice that $\phi(D)\vf=(\mathcal{F}^{-1}\phi) * \vf$, and for $t,u\in \bR^n$, we have the sampling representation (see, e.g., \cite[Section 3]{Frazier:2004ub} or \cite[Section 6]{MR1107300}), $$(\mathcal{F}^{-1}\phi * \vf)(t)=\sum_{\ell\in \bZ^n} \vf(\ell+u)[\mathcal{F}^{-1}\phi](t-u-\ell).$$
Hence,
$$W^{1/p}(t)(\mathcal{F}^{-1}\phi * \vf)(t)=\sum_{\ell\in \bZ^n} W^{1/p}(t)\vf(\ell+u)[\mathcal{F}^{-1}\phi](t-u-\ell),$$
which implies that,
\begin{align*}
|W^{1/p}(t)(\mathcal{F}^{-1}\phi * \vf)(t)|\leq \sum_{\ell\in \bZ^n} |W^{1/p}(t)\vf(\ell+u)[\mathcal{F}^{-1}\phi](t-u-\ell)|.
\end{align*}
As before, we let $Q_k:=Q(k,1)$ for $k\in\bZ^n$. Using $0<p\leq 1$, we arrive at the following estimate valid for $u\in Q_0$, $t\in\bR^n$,
\begin{align}
|W^{1/p}(t)(\mathcal{F}^{-1}\phi * \vf)(t)|^p&\leq \sum_{\ell\in \bZ^n} |W^{1/p}(t)\vf(\ell+u)|^p|[\mathcal{F}^{-1}\phi](t-u-\ell)|^p\nonumber\\
&\leq K^p\sum_{\ell\in \bZ^n} |W^{1/p}(t)\vf(\ell+u)|^p(1+|t-u-\ell|)^{-Mp}\nonumber\\
&\leq c_MK^p\sum_{\ell\in \bZ^n} |W^{1/p}(t)\vf(\ell+u)|^p(1+|t-\ell|)^{-Mp}.
\label{eq:norm}
\end{align}
 We now average the estimate \eqref{eq:norm} over $u\in Q_0$, noticing that $|Q_0|=1$,
\begin{align}
|W^{1/p}(t)(\mathcal{F}^{-1}\phi * \vf)(t)|^p&\leq c_MK^p\int_{Q_0}\sum_{\ell\in \bZ^n} 
|W^{1/p}(t)\vf(\ell+u)|^p(1+|t-\ell|)^{-Mp}\,du\nonumber\\
&= c_MK^p\sum_{\ell\in \bZ^n} (1+|t-\ell|)^{-Mp} \int_{Q_0}
|W^{1/p}(t)\vf(\ell+u)|^p\,du\nonumber\\
&= c_MK^p\sum_{\ell\in \bZ^n} (1+|t-\ell|)^{-Mp} \int_{Q_\ell}
|W^{1/p}(t)\vf(y)|^p\,dy,\label{eq:esti2}
\end{align}
where we have used Tonelli's theorem.  By the  doubling condition satisfied by the scalar weight $w_{\vw}(t):=|W^{1/t}(t)\vw|^p$, with doubling exponent $\beta$ independent of $\vw\in\bC^N$, there exists a finite constant $c_w$ such that for any $k,\ell\in\bZ^n$, and $y\in\bR^n$,
\begin{equation}\label{eq:doub}
\int_{Q_k} |W^{1/p}(t)\vf(y)|^p\,dt\leq c_w(1+|k-\ell|)^\beta\int_{Q_\ell} |W^{1/p}(t)\vf(y)|^p\,dt.
\end{equation}
 For $k\in\bZ^n$, we integrate inequality \eqref{eq:esti2} over $t\in Q_k$, using Tonelli's theorem once more together with the estimate \eqref{eq:doub}, and the observation that there exists $c>0$ such that for all $t\in Q_k$ and $\ell\in \bZ^n$,
$1+|k-\ell|\leq c(1+|t-\ell|)$,
\begin{align*}
\int_{Q_k} |W^{1/p}(t)(\phi(D)\vf)(t)|^p\,dt&\leq c_Mc^{Mp}K^p
 \int_{Q_k}\sum_{\ell\in \bZ^n} (1+|k-\ell|)^{-Mp} \int_{Q_\ell}
|W^{1/p}(t)\vf(y)|^p\,dy\,dt\\
&\leq c_M'\sum_{\ell\in \bZ^n} (1+|k-\ell|)^{-Mp+\beta} \int_{Q_\ell}\int_{Q_\ell}
|W^{1/p}(t)\vf(y)|^p\,dt\,dy,
\end{align*}
with $c_M':=c_wc_Mc^{Mp}K^p$. In the following concluding estimate, the matrix $A_p$ condition \eqref{eq:mA1} for $W$ will be essential. We first notice that $\{Q_k\}_k$ forms a partition of $\bR^n$ with $|Q_k|=1$, so we have, using the assumption that $Mp-\beta>n$ and putting 
$L:=\sum_k(1+|k|)^{-Mp+\beta}<\infty$,
\begin{align*}
\|\phi(D)\vf\|_{L^p(W)}^p&=\sum_{k\in \bZ^n} \int_{Q_k} |W^{1/p}(t)(\phi(D)\vf)(t)|^p\,dt\\
&\leq c_M'\sum_{k\in \bZ^n}\sum_{\ell\in \bZ^n} (1+|k-\ell|)^{-Mp+\beta}\int_{Q_\ell}\int_{Q_\ell}
|W^{1/p}(t)\vf(y)|^p\,dt\,dy\\
&=Lc_M'\sum_{\ell\in \bZ^n} \int_{Q_\ell}\int_{Q_\ell}
|W^{1/p}(t)\vf(y)|^p\,dt\,dy\\
&=Lc_M'\sum_{\ell\in \bZ^n} \int_{Q_\ell}\int_{Q_\ell}
|W^{1/p}(t)W^{-1/p}(y)W^{1/p}(y)\vf(y)|^p\,dt\,dy\\
&\leq Lc_M'\sum_{\ell\in \bZ^n}\int_{Q_\ell}\bigg( \int_{Q_\ell}
\|W^{1/p}(t)W^{-1/p}(y)\|^p\,dt\bigg) |W^{1/p}(y)\vf(y)|^p\,dy\\
&= Lc_M'\sum_{\ell\in \bZ^n}\int_{Q_\ell}\bigg( \frac{1}{|Q_\ell|}\int_{Q_\ell}
\|W^{1/p}(t)W^{-1/p}(y)\|^p\,dt\bigg) |W^{1/p}(y)\vf(y)|^p\,dy\\
&\leq Lc_M'  [W]_{\mathbf{A}_p(\bR^d)} \sum_{\ell\in \bZ^n}\int_{Q_\ell}|W^{1/p}(y)\vf(y)|^p\,dy\\
&=Lc_M'  [W]_{\mathbf{A}_p(\bR^d)}\|\vf\|_{L^p(W)}^p\\
&=C^p \|\vf\|_{L^p(W)}^p,
\end{align*}
with $C^p:=Lc_M'[W]_{\mathbf{A}_p(\bR^d)}$, which completes the proof in the case $R=1$. For  general $0<R<\infty$, we consider the multiplier $\psi:B(0,1)\rightarrow\bC$ defined by
$\psi(\cdot)=\phi(R\cdot)$. We  clearly have  
$$\mathcal{F}^{-1}(\psi)(x)=R^{-n}\mathcal{F}^{-1}(\phi)(R^{-1}x) \Longrightarrow 
|\mathcal{F}^{-1}(\psi)(x)|\leq K (1+|x|)^{-M}.$$
Now, notice that for any $\vf\in E_R$, we have $\vg:=R^{-n}\vf(R^{-1}\cdot)\in E_1$. Hence, using the result from the first part of the proof, we have for $\vf\in E_R$,
\begin{align*}
\|\phi(D)\vf\|_{L^p(W)}^p&=\int_{\bR^n} |W^{1/p}(t)\phi(D)\vf(t)|^p\,dt\\
&=R^{-n}\int_{\bR^n} |W^{1/p}(R^{-1}u)(\phi(D)\vf)(R^{-1}u)|^p\,du\\
&=R^{-n+np}\int_{\bR^n} |W^{1/p}(R^{-1}u)(\psi(D)\vg)(u)|^p\,du\\
&\leq C R^{-n+np}\int_{\bR^n} |W^{1/p}(R^{-1}u)\vg(u)|^p\,du\\
&=C R^{-n}\int_{\bR^n} |W^{1/p}(R^{-1}u)\vf(R^{-1}u)|^p\,du\\
&=C \int_{\bR^n} |W^{1/p}(t)\vf)(t)|^p\,dt,
\end{align*}
where $C:=C([W(R^{-1}\cdot)]_{\mathbf{A}_p(\bR^n)},K,p)$. However, as mentioned in Remark \ref{rem:inv},  the matrix ${A}_p$-condition is dilation invariant in the sense that $[W(R^{-1}\cdot)]_{\mathbf{A}_p(\bR^n)}=[W]_{\mathbf{A}_p(\bR^n)}$ for any $R>0$, making $C$ independent of $R$. This completes the proof. 
\end{proof}

\begin{remark}\label{rem:co}
As mentioned in the introduction, in case $1\leq p<\infty$,  a corresponding multiplier result holds without any assumptions on the support of the multiplier nor on the spectrum of $\vf$, see \cite[Lemma 4.4]{MR4263690}:   Suppose   $W\in\mathbf{A}_p(\bR^n)$ and assume there is a constant $K$ such that $\phi:\bR^n\rightarrow\bC$ satisfies 
\begin{equation}\label{eq:dec2}
|\mathcal{F}^{-1}(\phi)(x)|\leq KR^{n} (1+R|x|)^{-n-1},\qquad x\in\bR^n,
\end{equation}
then there exists a finite constant $C:=C([W]_{\mathbf{A}_p(\bR^n)},K,p)$ such that
$$\|\phi(D)\vf\|_{L^p(W)}\leq C \|\vf\|_{L^p(W)},\qquad \vf\in L^p(W).$$
\end{remark}

\section{An application: matrix-weighted smoothness spaces}\label{sec:exa}
An important application of   matrix-weighted $L^p$-spaces
is to the construction of various matrix-weighted smoothness spaces obtained by  imposing suitable weighted $L^p$-restrictions on local components of a vector-function. The  local components  are often defined with an aim  to capture local frequency content of the (vector-)function and  can naturally be  obtained by applying a family of suitable bandlimited Fourier multipliers (frequency filters) compatible with a desired decomposition of the frequency space.

Roudenko was the first to apply such an approach in the matrix weighted setup, see \cite{Rou03a}, where she introduced a very natural notion of matrix-weighted Besov spaces $B^s_{p,q}(W)$ based on dyadic decompositions of vector-functions. This work was later extended by Frazier and Roudenko  \cite{Frazier:2004ub, MR4263690} to matrix-weighted Triebel-Lizorkin spaces and Besov spaces for the full range $0<p<\infty$. A very extensive recent study by Bu et al.\ of various properties of matrix-weighted Besov and Triebel-Lizorkin type spaces can be found in \cite{bu2023matrixweighted,bu2023matrixweighted2,MR4797823}.

We shall only touch upon one aspect of this rather involved theory here, namely the issue of such smoothness spaces being well-defined.

Let $W\in \mathbf{A}_p(\bR^n)$ for some $0<p<\infty$, and let $\beta>0$ be the doubling exponent from Eq.\ \eqref{eq:dou} associated with $W$. 

We say that $\Psi:=\{\psi_j\}_{j\in\bZ}$ is a bounded  admissible partition of unity if it is a smooth resolution of the identity on $\bR^n\backslash \{0\}$, i.e., $\sum_j \psi_j(x)\equiv 1$ on $\bR^n\backslash \{0\}$, and there exist $c_1:=c_1(\Psi), c_2:=c_2(\Psi)>0$ such that $\text{supp}(\psi_j)\subseteq \{x\in\bR^n:c_12^j\leq |x|<c_2 2^j\}$, $j\in\bZ$. Also suppose there exists $C:=C(\Psi)$ such that
\begin{equation}\label{eq:dec}|\mathcal{F}^{-1}\psi_j(x)|\leq C2^{jn}(1+2^{j}|x|)^{-M},\qquad j\in\bZ,
\end{equation}
for some $M>(n+\beta)/\min\{1,p\}$. Then we follow Frazier and Roudenko and define the homogeneous
matrix-weighted Besov space $\dot{B}_{p,q}^s(W):=\dot{B}_{p,q}^s(W,\Psi)$ for $0<q	\leq \infty$ as the collection of $\vf=(f_1,\ldots f_N)$ with $f_i\in \ddS/\mathcal{P}$ (the tempered distributions modulo polynomials), $i=1,\ldots,N$, satisfying
$$\|\vf\|_{\dot{B}_{p,q}^s(W)}:=\bigg(\sum_{j=-\infty}^\infty 2^{jsq}\|\psi_j(D)\vf\|_{L^p(W)}^q\bigg)^{1/q}<+\infty,$$
with the sum replaced by $\sup_j$ in the case $q=\infty$.
Now, let us take another bounded admissible partition of unity $\Phi:=\{\phi_j\}_{j\in\bZ}$. Clearly, it is desirable that the construction of $\dot{B}_{p,q}^s(W)$ is well-defined in the sense  that $\dot{B}_{p,q}^s(W,\Psi)=\dot{B}_{p,q}^s(W,\Phi)$ up to equivalence of norms, and this important fact is indeed proved in \cite{Frazier:2004ub,Rou03a}, but the reader may verify that the proofs presented in \cite{Frazier:2004ub,Rou03a} of this fact are somewhat involved relying on a  theory of almost diagonal matrices in the matrix-weighted setting developed by the same authors. We now give an alternative, more transparent, proof of the mentioned equivalence relying {\em solely} on Proposition \ref{prop:main}.

We first notice that for $j\in \bZ$,
\begin{equation}\label{eq:pdo}
\psi_j(D)\vf=\psi_j(D)\sum_{k\in A_j}\phi_k(D)\vf,
\end{equation}
with $A_j=\{k\in \bZ:\text{supp}(\psi_j)\cap \text{supp}(\phi_k)\not=\emptyset\}$, where it is easy to verify that $\#A_j$ is bounded by a constant $n_0$ independent of $j$ due to the dyadic nature of the support sets in $\Psi$ and $\Phi$. Hence, using the decay assumption \eqref{eq:dec}, we may call on Proposition \ref{prop:main} in case $0<p\leq 1$. In case $1<p<\infty$, we may use the result mentioned in Remark~\ref{rem:co}, where we notice that $M>n+\beta\geq n+1$ since $\beta\geq n\geq 1$, so the condition stated in Eq.~\eqref{eq:dec2} is satisfied. For any $0<p<\infty$, we obtain
$$\|\psi_j(D)\vf\|_{L^p(W)}\leq C\sum_{k\in A_j}\|\phi_k(D)\vf\|_{L^p(W)}.$$
Also observe that for $k\in A_j$, we have $2^{ks}\asymp 2^{js}$ uniformly in $j$.
 It follows from this observation that
$$2^{js}\|\psi_j(D)\vf\|_{L^p(W)}\leq C\sum_{k\in A_j} 2^{ks}\|\phi_k(D)\vf\|_{L^p(W)}.$$
Using the uniform bounds on the cardinality of the sets $A_j$, it is  then straightforward to verify that
\begin{align*}
\|\vf\|_{\dot{B}^{s}_{p,q}(W,\Psi)}&:=\bigg(\sum_{j=-\infty}^\infty 2^{jsq}\|\psi_j(D)\vf\|_{L^p(W)}^q\bigg)^{1/q}\\
&\leq C \bigg(\sum_{j=-\infty}^\infty \bigg(\sum_{k\in A_j} 2^{ks}\|\phi_k(D)\vf\|_{L^p(W)}\bigg)^q\bigg)^{1/q}\\
&\leq C' \bigg(\sum_{j=-\infty}^\infty \sum_{k\in A_j} 2^{ksq}\|\phi_k(D)\vf\|_{L^p(W)}^q\bigg)^{1/q}\\
&\leq C'' \bigg(\sum_{k=-\infty}^\infty 2^{ksq}\|\phi_k(D)\vf\|_{L^p(W)}^q\bigg)^{1/q}\\
&=:C''\|\vf\|_{\dot{B}^{s}_{p,q}(W,\Phi)}.
\end{align*}
Interchanging the roles of $\Psi$ and $\Phi$ provides the reverse estimate, yielding the wanted norm equivalence between $\dot{B}_{p,q}^s(W,\Psi)$ and $\dot{B}_{p,q}^s(W,\Phi)$.

\end{document}